\input amstex
\documentstyle{amsppt}

\topmatter
\title Integrability of almost complex structures on Banach manifolds\endtitle
\author Daniel Belti\c t\u a\endauthor
\address Institute of Mathematics ``Simion Stoilow'' of the Romanian Academy, 
P.O. Box 1-764, RO-014700 Bucharest, Romania\endaddress
\email Daniel.Beltita\@ imar.ro\endemail
%\date July 14, 2004\enddate
%\keywords almost complex structure; Banach-Lie group; homogeneous space\endkeywords
%\subjclass Primary 32Q60; Secondary 53C15;58B12\endsubjclass   
\abstract We prove that the classical integrability condition for 
almost complex structures 
on finite-dimensional smooth manifolds also works in infinite dimensions 
in the case of almost complex structures that are real analytic 
on real analytic Banach manifolds. 
As an application, we extend some known results concerning 
existence of invariant complex structures on homogeneous spaces of Banach-Lie groups. 
 

\endabstract
%\leftheadtext{}
\rightheadtext{Almost complex structures on Banach manifolds}
\endtopmatter

\def\1{\text{\bf 1}}
\def\Ad{\text{\rm Ad}}
\def\Adt{\widetilde{\text{\rm Ad}}}
\def\Aut{\text{\rm Aut}}
\def\Fl{\text{\rm Fl}}
\def\gg{{\goth g}}
\def\GL{\text{\rm GL}}
\def\hg{{\goth h}}
\def\kg{{\goth k}}
\def\mult{\text{\rm mult(i)}}
\def\Ic{{\Cal I}}
\def\Iac{{\Cal I}\text{\it ac}}
\def\Ico{{\Cal I}\text{\it co}}
\def\Bc{{\Cal B}}
\def\Kc{{\Cal K}}
\def\Oc{{\Cal O}}
\def\Vc{{\Cal V}}
\def\Zc{{\Cal Z}}
\def\Ca{{\Cal C}^\omega}
\def\Ci{{\Cal C}^\infty}
\def\Ker{\text{\rm Ker}\,}
\def\L{\text{\bf L}}
\def\Ran{\text{\rm Ran}\,}
\def\er{\text{\rm e}}
\def\id{\text{\rm id}}
\def\ir{\text{\rm i}}
\def\U{\text{\rm U}}

\document

\head 1. Introduction\endhead

In this paper we prove that the classical integrability condition for 
almost complex structures 
on finite-dimensional smooth manifolds (see \cite{EF51}, \cite{Fr55}, \cite{NN57})
also works in infinite dimensions, ensuring the existence of complex coordinates 
in the case of almost complex structures that are real analytic on 
real analytic Banach manifolds. 
And this might turn out to be interesting in connection with both 
the recent interest in fundamental problems in infinite-dimensional complex analysis 
(see \cite{Le98}, \cite{Le99} and \cite{Le00}) 
and the role of infinite-dimensional complex homogeneous spaces 
in representation theory of infinite-dimensional Lie groups 
(see for instance \cite{NRW01} and \cite{Ne04}).

Our motivation in writing this paper came from the fact that, 
in the existing literature, 
the subject of
complex coordinates in almost complex Banach manifolds was approached in two very different manners: 
either one constructed the complex coordinates in a way that seemed to be very specific to the manifolds under consideration
(see for instance 
\cite{Ne04}, \cite{BR04}, \cite{Be03}, \cite{MS98}, \cite{MS97}, \cite{Wi94} and 
\cite{Wi90}) 
or one exhibited examples of almost complex structures that do not come from any complex coordinates  at all 
(see \cite{Pa00}). 
Our point is to show that there exists a large class of almost complex structures, 
namely the torsion-free real analytic almost complex structures on real analytic Banach manifolds, that always come from complex coordinates on the corresponding manifold 
(Theorem 7 below). 
That this class of almost complex structures is really wide can be seen from the fact that all 
invariant almost complex structures on homogeneous spaces of Banach-Lie groups are real analytic (see Proposition~11 below). 

As a matter of fact, our main result on invariant complex structures (Theorem~15) works for homogeneous spaces of not necessarily connected Banach-Lie groups. 
This point is particularly important for 
the applications to the adjoint or coadjoint orbits that show up in connection with 
the operator algebras 
(see \cite{BR04}, \cite{MS98}, \cite{MS97} and Corollary~16 below), 
since the unitary groups of $C^*$-algebras are in general non-connected 
(although they are connected in the special case of $W^*$-algebras), 
while nothing seems to be known on connectedness of the unitary groups associated with 
non-separable operator ideals 
(compare Remark~4.9 in~\cite{BR04}). 
Thus, our method improves the one suggested in \cite{Ne04}, 
where one needs the connectedness of the group that transitively acts on the homogeneous space.
(The method suggested in \cite{Ne04} is based on extending to infinite dimensions certain reasonings from \cite{Ki76} 
involving local groups and Baker-Campbell-Hausdorff series.)

The basic idea in the proofs of our main results is to make sure that the main steps in the finite-dimensional original reasonings (see \cite{Fr55}) carry over to infinite dimensions. 
In order to do so, a few technical devices are needed from infinite-dimensional global analysis (see \cite{Up85}, \cite{Nu92} and \cite{La01}). 
The paper is organized as follows: 
In Section~2 we prove the main result on almost complex structures on general manifolds 
(Theorem~7). 
In Section~3 we consider the case of homogeneous spaces of Banach-Lie groups, 
and the main result of this section 
describes a way to parameterize the invariant complex structures in Lie algebraic terms 
(Theorem~15).

We now introduce some notation that will be used throughout the paper. 
We denote by $\ir$ the square root of $-1$, 
by $\er$ the basis of the natural logarithms, 
by $TM$ the tangent bundle of a manifold $M$
and by  $B_X(a,r)$ the open ball of center $a\in X$ and radius $r>0$ 
in the real Banach space $X$. 
Moreover $\Bc(X)$ stands for the bounded real-linear operators on the real Banach 
space $X$, while  $\Bc^{\Bbb C}(Z)$ denotes the bounded complex-linear operators on a complex Banach space $Z$. 
As regards function spaces, 
$\Ci(D,Y)$ is the set of $\Ci$ functions on the open subset $D$ of some real Banach space with values in another real Banach space $Y$, 
 $\Ci(D):=\Ci(D,{\Bbb R})$
$\Ca(D,Y)$ is the space of real analytic functions on the open subset $D$ of some real Banach space with values in the real Banach space $Y$,
$\Oc(D,W)$ is the space of holomorphic functions on the open subset $D$ of some complex Banach space with values in 
another complex Banach space $W$, and 
$\Oc(D):=\Oc(D,{\Bbb C})$. 

We denote by $X_{\Bbb C}=X\dotplus\ir X$
the complexified Banach space of the real Banach space $X$. 
We always think of $X_{\Bbb C}$ equipped with a norm such that 
$\|x-\ir y\|=\|x+\ir y\|$ for all $x,y\in X$. 
Moreover, unless otherwise stated, we denote by $z\mapsto\overline{z}$ the 
conjugation of $Z$ whose set of fixed points equals $X$, that is, 
$\overline{x+\ir y}=x-\ir y$ whenever $x,y\in X$. 
We also note that every $T\in\Bc(X)$ can be uniquely extended to an element 
of $\Bc^{\Bbb C}(X_{\Bbb C})$ which will always be denoted by $T$ as well. 

For an open subset $D$ of a real Banach space $X$,  
we consider the usual Lie bracket $[\cdot,\cdot]$ defined on $\Ci(D,X)$ 
by 
$$[a,b]_x=a'_x(b_x)-b'_x(a_x)$$
whenever $x\in D$ and $a,b\in\Ci(D,X)$, where we denote the value of $a$ at $x$ by 
$a_x:=a(x)\in X$, the differential of $a$ at $x$ by $a'_x\colon X\to X$, 
the second differential of $a$ at $x$ by $a''_x\colon X\times X\to X$, and so on. 
Also, for a smooth mapping $J\colon D\to\Bc(X)$ and any $a\in\Ci(D,X)$ 
we define  $Ja\in\Ci(D,X)$ by $(Ja)_x=J_xa_x$ for all $x\in D$.

\head 2. General results\endhead

\definition{Definition 1}
Let $M$ be a smooth real Banach manifold. 
An {\it almost complex structure} on $M$ is a mapping 
$J\colon TM\to TM$ such that 
for all $p\in M$ we have $J(T_pM)\subseteq T_pM$, 
$J_p:=J|_{T_pM}\colon T_pM\to T_pM$ is a bounded (real-)linear operator 
and $(J_p)^2=-\id_{T_pM}$.
\enddefinition

\proclaim{Lemma 2} 
Let $X$ be a real Banach space and $D$ an open subset of $X$. 
Suppose that 
$$J\colon D\to\Bc(X),\quad x\mapsto J_x$$
is a smooth mapping satisfying $(J_x)^2=-\id_X$ for all $x\in D$, and 
consider 
$$\Omega\colon\Ci(D,X)\times\Ci(D,X)\to\Ci(D,X)$$
defined by 
$$\Omega(a,b)=J[Ja,b]+J[a,Jb]+[a,b]-[Ja,Jb]$$
for all $a,b\in\Ci(D,X)$. 
Then we have 
$$\Omega(a,b)_x=-J'_x(J_xa_x,b_x)-J'_x(a_x,J_xb_x)+J'_x(J_xb_x,a_x)+J'_x(b_x,J_xa_x)$$
for all $a,b\in\Ci(D,X)$ and $x\in D$. 
\endproclaim

\demo{Proof}
The proof is a straightforward computation based on the formulas 
$$(Jc)'_xu=J'_x(c_x,u)+J_xc'_xu
\quad\text{whenever }x\in D,u\in X\text{ and }c\in\Ci(D,X),
\eqno{(1)}$$
and 
$$J'_x(J_xu,v)+J_xJ'_x(u,v)=0\quad\text{whenever }x\in D\text{ and }u,v\in X.
\eqno{(2)}$$
We omit the details. 
\qed
\enddemo

\definition{Definition 3}
Let $M$ be a smooth real Banach manifold,  
denote by $TM\oplus TM$ the Whitney sum of the vector bundle $TM$ with itself 
and assume that $J\colon TM\to TM$ is a smooth almost complex structure. 
Then the {\it torsion} of $J$ is the smooth mapping 
$$\Omega\colon TM\oplus TM\to TM$$
having the property that, if $a$ and $b$ are locally defined smooth vector fields on $M$, 
then we have 
$$\Omega(a,b)=J[Ja,b]+J[a,Jb]+[a,b]-[Ja,Jb]$$
on any open set where both $a$ and $b$ are defined. 
Note that the existence of the mapping $\Omega$ is a consequence of Lemma~2. 
\enddefinition

\proclaim{Lemma 4}
Let $X$ be a real Banach space, $r>0$ and $J\colon B_X(0,r)\to\Bc(X)$ 
a real analytic mapping satisfying $(J_x)^2=-\id_X$ whenever $x\in B_X(0,r)$. 
Consider the complexified Banach space $Z=X_{\Bbb C}$ and extend 
$J|_{B_X(0,r/2\er)}$ to a holomorphic mapping 
$J\colon B_Z(0,r/2\er)\to\Bc^{\Bbb C}(Z)$. 
Next denote 
$$\Zc_{\pm}:=\{a\in\Oc(B_Z(0,r/2\er),Z)\mid Ja=\mp\ir a\}$$
and define 
$$\Omega\colon\Ci(D,Y)\times\Ci(D,Y)\to\Ci(D,Y)$$
as in {\rm Lemma~2}, where $Y$ stands for either of the spaces $X$ and $Z$, 
and $D$ is any open subset of the corresponding Banach space
($Z$ is thought of as a real Banach space). 

Now  consider the following assertions: 

\itemitem{$1^\circ$} For all $a,b\in\Ca(B_X(0,r/2\er),X)$ we have $\Omega(a,b)=0$. 

\itemitem{$2^\circ$} For all $a,b\in\Zc_{-}$ we have $[a,b]\in\Zc_{-}$. 

\itemitem{$3^\circ$} For all $a,b\in\Zc_{+}$ we have $[a,b]\in\Zc_{+}$.

\itemitem{$4^\circ$} For all $a,b\in\Ca(B_X(0,r),X)$ we have $\Omega(a,b)=0$. 

\noindent Then we have 
$1^\circ\Rightarrow 2^\circ \Leftrightarrow 3^\circ \Rightarrow 2^\circ$.
\endproclaim

\demo{Proof}
We first note that, if $n\in{\Bbb Z}_{+}$, $f_n\colon Z\times\cdots\times Z\to Z$ 
is complex $n$-linear and $f_n(X\times\cdots\times X)\subseteq X$, 
then 
$$\aligned
\|f_n(x_1+\ir y_1,\dots,x_n+\ir y_n)\|
 &\le \|f_n|_{X\times\cdots\times X}\|\cdot(\|x_1\|+\|y_1\|)\cdots(\|x_n\|+\|y_n\|) \cr
 &\le 2^n\|f_n|_{X\times\cdots\times X}\|\cdot(\|x_1+\ir y_1\|)\cdots(\|x_n+\ir y_n\|),
\endaligned$$
hence $\|f_n\|\le 2^n\|f_n|_{X\times\cdots\times X}\|$. 
It follows by this remark along with Corollary~1.5 in \cite{Up85} 
that for each $a\in\Ca(B_X(0,r),X)$ there is a mapping in $\Oc(B_Z(0,r/2\er),Z)$ 
that agrees with $a$ on $B_X(0,r/2\er)$; that holomorphic mapping 
is uniquely determined according to Theorem~1.11 in \cite{Up85}, 
and therefore it will be denoted by $a$ as well. 
Similarly we define $J\colon B_Z(0,r/2\er)\to\Bc^{\Bbb C}(Z)$. 
Throughout the proof, we also define 
$$\Omega\colon \Oc(B_Z(0,r/2\er),Z)\times \Oc(B_Z(0,r/2\er),Z)
\to\Oc(B_Z(0,r/2\er),Z)$$ 
by the same formula as in Lemma~2. 
Now we come back to the proof of the asserted relationship between assertions 
$1^\circ$--$4^\circ$. 

$1^\circ\Rightarrow 2^\circ$ 
Let $a,b\in\Zc_{-}$, so that $Ja=-\ir a$ and $Jb=-\ir b$. 
By $1^\circ$ we have $\Omega(a,b)|_{B_X(0,r/2\er)}=0$, 
hence $\Omega(a,b)=0$ on $B_Z(0,r/2\er)$ by Theorem~1.11 in~\cite{Up85}. 
Then 
$$0=J[-\ir a,b]+J[a,-\ir b]+[a,b]-[-\ir a,-\ir b]=-2\ir J[a,b]+2[a,b],$$
whence $J[a,b]=-\ir[a,b]$ as desired. 

$2^\circ\Rightarrow 3^\circ$ 
If $a,b\in\Zc_{-}$, then clearly $\overline{a(\cdot)},\overline{b(\cdot)}\in\Zc_{+}$, 
hence $[\overline{a(\cdot)},\overline{b(\cdot)}]\in\Zc_{+}$. 
Now $\overline{[a,b](\cdot)}=[\overline{a(\cdot)},\overline{b(\cdot)}]\in\Zc_{+}$, 
whence $[a,b]\in\Zc_{-}$. 

$3^\circ\Rightarrow 2^\circ$ 
Similar to $2^\circ\Rightarrow 3^\circ$. 

$3^\circ\Rightarrow 4^\circ$ 
Let $a,b\in\Ca(B_X(0,r),X)$ arbitrary and construct the corresponding 
$a,b\in\Oc(B_Z(0,r/2\er),Z)$ as explained at the beginning of the proof. 
Then we have 
$$a_{\pm}:=\frac{1}{2}(a\mp\ir Ja)\in\Zc_{\pm}, \quad
  b_{\pm}:=\frac{1}{2}(b\mp\ir Jb)\in\Zc_{\pm},$$
and 
$$a=a_{+}+a_{-}\text{ and }b=b_{+}+b_{-}\text{ on }B_Z(0,r/2\er).$$
In particular, 
$$\Omega(a,b)= \Omega(a_{+},b_{+})
              +\Omega(a_{+},b_{-})
              +\Omega(a_{-},b_{+})
              +\Omega(a_{-},b_{-})
\text{ on }B_Z(0,r/2\er).$$
Now, since $Ja_{\pm}=\mp\ir a_{\pm}$ and $Jb_{\pm}=\mp\ir b_{\pm}$, 
we easily get by the definition of $\Omega(\cdot,\cdot)$ that 
$\Omega(a,b)=0$ on $B_Z(0,r/2\er)$. 
Then Theorem~1.11 in \cite{Up85} 
shows that $\Omega(a,b)=0$ on $B_X(0,r)$, 
and the proof is finished. 
\qed
\enddemo

\proclaim{Lemma 5}
In the setting of {\rm Lemma~4}, denote 
$$Z_0^{-}:=\Ker(J_0+\ir\cdot\id_Z)\quad(\subseteq Z).$$
Then there exists $r_1\in(0,r/2\er)$ such that for every $z\in B_Z(0,r_1)$ 
the sequence 
$$0 @>>> Z_0^{-} @>{J_z-\ir\cdot\id_Z}>> Z @>{J_z+\ir\cdot\id_Z}>> Z$$
is exact.
\endproclaim

\demo{Proof}
To begin with, denote 
$Z_z^{-}:=\Ker(J_z+\ir\cdot\id_Z)\subseteq Z$ 
for all $z\in B_Z(0,r/2\er)$. 
Since $(J_z)^2=-\id_Z$, it then follows that 
for all $z\in B_Z(0,r/2\er)$ we have 
$$Z_z^{-}=\Ran(J_z-\ir\cdot\id_Z)\text{ and }
Z=Z_z^{-}\dotplus\Ker(J_z-\ir\cdot\id_Z).$$
In particular, it follows that 
$(J_0-\ir\cdot\id_Z)|_{Z_0^{-}}\colon Z_0^{-}\to Z_0^{-}$
is an isomorphism of the complex Banach space $Z_0^{-}$ onto itself. 

On the other hand, for arbitrary $z\in B_Z(0,r)$, 
note that we have an isomorphism of real Banach spaces 
$$\Psi_z\colon Z_z^{-}\to X,\quad z\mapsto\frac{1}{2}(z+\overline{z}).$$
In fact, since $\overline{J_zv}=J_z\overline{v}$ for all $v\in Z$, 
it is easy to check that $\Psi_z$ 
has an inverse defined by $x\mapsto x+\ir J_z x$. 

Now we have a continuous mapping 
$$\tau\colon B_Z(0,r)\to\Bc(Z_0^{-},X),\quad 
z\mapsto\Psi_z\circ (J_z-\ir\cdot\id_Z)|_{Z_0^{-}},$$
and this mapping has 
the property that $\tau(0)\colon Z_0^{-}\to X$ is an isomorphism, 
since both 
$(J_0-\ir\cdot\id_Z)|_{Z_0^{-}}\colon Z_0^{-}\to Z_0^{-}$ 
and 
$\Psi_0\colon Z_0^{-}\to X$ are isomorphisms. 
Since the isomorphisms $Z_0^{-}\to X$ constitute 
an open subset of $\Bc(Z_0^{-},X)$, 
it then follows that there exists $r_1\in(0,r/2\er)$ such that 
$\tau(z)\colon Z_0^{-}\to X$ is an isomorphism of real Banach spaces 
whenever $z\in B_Z(0,r_1)$. 
Since each $\Psi_z\colon Z_z^{-}\to X$ is an isomorphism, 
we deduce that 
$$(J_z-\ir\cdot\id_Z)|_{Z_0^{-}}\colon Z_0^{-}\to Z_z^{-}=\Ker(J_z+\ir\cdot\id_Z)$$
is an isomorphism for all $z\in B_Z(0,r_1)$, 
and the proof ends. 
\qed
\enddemo

\proclaim{Theorem 6}
Let $M$be a complex Banach manifold and $E$ a complex subbundle of $TM$ 
with the following property: if $V$ is an open subset of $M$ and 
$\xi,\eta\colon V\to TM$ are holomorphic vector fields on $V$ such that 
$\xi(p),\eta(p)\in E_p$ for all $p\in V$, 
then $[\xi,\eta](p)\in E_p$ for all $p\in V$. 
Then for every $m\in M$ there exist an open neighborhood $U$ of $m$, 
a complex Banach space $F$ and a holomorphic submersion $\zeta\colon U\to F$ 
such that 
$$E_p=\Ker(T_p\zeta)=T_p(\zeta^{-1}(\zeta(p)))\text{ for all }p\in U.$$
\endproclaim

\demo{Proof}
All the needed reasonings that are carried out in Chapter~VI in \cite{La01} 
in order to prove the Frobenius theorem for ${\Cal C}^k$ mappings 
on real Banach spaces work for holomorphic mappings on complex Banach spaces as well. 
\qed
\enddemo

Now we are ready to prove the following infinite-dimensional 
version of the result of \S 11 in~\cite{Fr55}. 
The proof is inspired by the last remark in \S 16 in~\cite{Fr55}. 

\proclaim{Theorem 7}
Let $C$ be a real analytic Banach manifold and 
$$J\colon TC\to TC$$ 
a real analytic almost complex structure on $C$. 
Denote by 
$$\Omega\colon TC\oplus TC\to TC$$
the torsion of $J$. 
Then $\Omega=0$ if and only if there exists on $C$ a structure of complex 
Banach manifold whose underlying real analytic structure is 
just the original structure of $C$ and such that 
$J_pv=\ir v$ whenever $p\in C$ and $v\in T_pC$. 
\endproclaim

\demo{Proof}
The proof has several stages. 

$1^\circ$
It is clear that $\Omega=0$ is a necessary condition for 
the existence of a complex structure as indicated in the statement. 

$2^\circ$ 
Conversely, assume that $\Omega=0$. 
Since the statement has a local character, 
it suffices to prove the following fact. 
Let $X$ be a real Banach space, $r>0$ and 
$$J\colon B_X(0,r)\to\Bc(X),\quad x\mapsto J_x,$$
a real analytic mapping such that $(J_x)^2=-\id_X$ for all $x\in B_X(0,r)$ 
and 
$$J[a,Jb]+J[Ja,b]+[a,b]-[Ja,Jb]=0 
\quad\text{for all }a,b\in\Ca(B_X(0,r),X).$$
Then there exist $r_0\in(0,r)$, a complex manifold $N$ and 
a real analytic diffeomorphism $\Theta\colon B_X(0,r_0)\to N$ 
such that 
$$T_x\Theta\circ J_x=\ir T_x\Theta
\quad\text{for all }x\in B_X(0,r_0).$$

$3^\circ$ 
To prove the assertion from stage~$2^\circ$, let 
$Z$, $J\colon B_Z(0,r/2\er)\to\Bc^{\Bbb C}(Z)$ and $\Zc_{\pm}$ as in Lemma~4. 
Moreover denote $M=B_Z(0,r/2\er)$ and 
$$E_z=\{z\}\times\Ker(J_z+\ir\cdot\id_Z)\subseteq M\times Z\quad
\text{for all }z\in M.$$
It then follows by Lemma~5 that $E=\bigcup\limits_{z\in M}E_z$ 
is a complex subbundle of $TM$ 
(see \cite{La01}). 

On the other hand, Lemma~4 ($1^\circ\Rightarrow3^\circ$) 
shows that, 
if $V$ is an open subset of $M$ and 
$\xi,\eta\colon V\to TM$ are holomorphic vector fields 
with $\xi(p),\eta(p)\in E_p$ for all $p\in V$. 
then $[\xi,\eta](p)\in E_p$ whenever $p\in V$. 
Now it follows by Theorem~6 that there exist $r_0\in(0,r/2\er)$, 
a complex Banach space $F$ and a holomorphic submersion 
$\zeta\colon B_Z(0,r_0)\to F$ 
such that 
$$\Ker(T_z\zeta)=E_z=\{z\}\times\Ker(J_z+\ir\cdot\id_Z)\quad
\text{for all }z\in B_Z(0,r_0).$$
In particular, since we have $\Ran(J_z-\ir\cdot\id_Z)=\Ker(J_z+\ir\cdot\id_Z)$, 
it follows that 
$T_z\zeta\circ(J_z-\ir\cdot\id_Z)=0$, 
whence 
$$T_z\zeta\circ J_z=\ir T_z\zeta\quad 
\text{ for all }z\in B_Z(0,r_0).$$

$4^\circ$ 
We now show that 
$$\Ker(J_0+\ir\cdot\id_Z)\dotplus X=Z.$$
In fact, since $J_0X=X$, it is clear that $\Ker(J_0+\ir\cdot\id_Z)\cap X=\{0\}$. 
Now let $z\in Z$ arbitrary. 
Then $z_{\pm}:=(J_0\pm\ir\cdot\id_Z)z/2\ir\in\Ker(J_0\mp\ir\cdot\id_Z)$ 
and $z=z_{+}-z_{-}$. 
Thus $z=(-z_{-}-\overline{z_{+}})+(z_{+}+\overline{z_{+}})$ 
with $-z_{-}-\overline{z_{+}}\in\Ker(J_0+\ir\cdot\id_Z)$ 
and 
$z_{+}+\overline{z_{+}}\in X$. 

Since $\Ker(T_0\zeta)=\{0\}\times\Ker(J_0+\ir\cdot\id_Z)$ 
and $\{0\}\times X=T_0(B_X(0,r_0))$, 
it then follows that 
$T_0(\zeta|_{B_X(0,r_0)})=(T_0\zeta)|_X\colon X\to\Ran(T_0\zeta)=F$ 
is an isomorphism of real Banach spaces. 
(Recall that $\zeta$ is a submersion.) 
The inverse mapping theorem then shows that 
there exist $r_1\in(0,r_0)$ and an open neighborhood $N$ 
of $\zeta(0)\in F$ such that 
$\Theta:=\zeta|_{B_X(0,r_1)}\colon B_X(0,r_1)\to N$ is a real analytic diffeomorphism. 
Now the last assertion in stage~$3^\circ$ shows that the diffeomorphism $\Theta$ 
has all the properties claimed in stage~$2^\circ$ of the proof. 
\qed
\enddemo

\head 3. The case of homogeneous spaces\endhead

\definition{Notation 8}
Throughout this section we denote by $G$ a real Banach-Lie group $G$ and 
by $H$ a Banach-Lie subgroup of $G$ (see for instance \cite{Up85}). 
Then $G/H$ has a structure of real analytic manifold such that the natural action 
$$\alpha\colon G\times G/H\to G/H,\quad (g,p)\mapsto\alpha(g,p)=\alpha_g(p)=g.p$$
is real analytic. 
Furthermore, we consider the natural projection 
$$\pi\colon G\to G/H,$$
which is a principal bundle. 
We denote the base point of $G/H$ by 
$$p_0:=H\in G/H.$$
The Lie algebras of $G$ and $H$ will be denoted by 
$\L(G)=\gg$ and $\L(H)=\hg$ respectively. 
Then we have a natural identification 
$$T_{p_0}(G/H)\simeq\gg/\hg.$$
Also note that the restriction to $H$ of the adjoint representation 
$\Ad_G\colon G\to\Aut(\gg)$ induces a continuous homomorphism of Banach-Lie groups 
$$\Adt_G\colon H\to\GL(\gg/\hg)=\GL(T_{p_0}(G/H)).$$
\enddefinition

\definition{Definition 9}
Let $J\colon T(G/H)\to T(G/H)$ be an almost complex structure 
on the homogeneous space $G/H$. 
We say that $J$ is {\it invariant} 
if for all $g\in G$ and $p\in G/H$ the diagram 
$$\CD
T_p(G/H) @>{T_p(\alpha_g)}>> T_{g.p}(G/H) \cr
@V{J_p}VV @VV{J_{g.p}}V \cr
T_p(G/H) @>{T_p(\alpha_g)}>> T_{g.p}(G/H) \cr
\endCD$$
is commutative.
\enddefinition

\proclaim{Proposition 10}
Consider the sets
$$\Iac(G,H)=\{J\colon T(G/H)\to T(G/H)\mid 
J\text{ invariant almost complex structure}\}$$
and 
$$\Ic_0(G,H)=\{I\in\Bc(\gg/\hg)\mid I^2=-\id_{\gg/\hg};\quad 
(\forall h\in H)\;\; T_{p_0}(\alpha_h)\circ I=I\circ T_{p_0}(\alpha_h)\}.$$
Then the mappings 
$$\check{\phantom{J}}\colon\Iac(G,H)\to\Ic_0(G,H),\quad J\mapsto \check{J}:=J_{p_0},$$
and 
$$\widehat{\phantom{I}}\colon\Ic_0(G,H)\to\Iac(G,H), \quad
\widehat{I}_{g.p_0}:=T_{p_0}(\alpha_g)\circ I\circ T_{p_0}(\alpha_g)^{-1} 
\text{ for all }g\in G,$$
are bijections inverse to one another.
\endproclaim

\demo{Proof}
For each $I\in\Ic_0(G,H)$ the mapping
$\widehat{I}_{g.p_0}\colon T_{g.p_0}(G/H)\to T_{g.p_0}(G/H)$ 
is well defined for all $g\in G$. 
Thus the mapping $\widehat{\phantom{I}}$ is well defined. 
Now the assertions are obvious. 
\qed
\enddemo

\proclaim{Proposition 11}
Every invariant almost complex structure $J$ on $G/H$ 
is real analytic. 
\endproclaim

\demo{Proof}
This follows by an easy application of Corollary~8.4 in \cite{Up85}. 
\qed
\enddemo

The following result is a more precise version of Satz~2 of \S 18 
in~\cite{Fr55} in infinite dimensions. 

\proclaim{Proposition 12}
Let $\Vc$ be a closed linear subspace of $\gg$ such that $\gg=\hg\dotplus\Vc$ 
and
$\Ic_{\Vc}(G,H)$ the set of all bounded linear operators 
$I\colon\gg\to\gg$ such that $\hg\subseteq\Ker I$, 
$I(\Vc)\subseteq\Vc$, $(I|_{\Vc})^2=-\id_{\Vc}$ 
and $T_{\1}\pi\circ I\circ\Ad_G(h)|_{\Vc}=T_{\1}\pi\circ\Ad_G(h)\circ I|_{\Vc}$ 
for all $h\in H$. 
Now for every $I\in\Ic_{\Vc}(G,H)$ define 
$c_{\Vc}(I)\colon\gg/\hg\to\gg/\hg$ by $c_{\Vc}(I)(x+\hg)=Ix+\hg$ 
for all $x\in\gg$. 

Then we obtain a bijection 
$$c_{\Vc}\colon\Ic_{\Vc}(G,H)\to\Ic_0(G,H).$$
\endproclaim

\demo{Proof}
The proof follows the lines of the proof of Satz~2 of \S 18 
in~\cite{Fr55}, 
so that we omit the details. 
\qed
\enddemo

We now prove an extension of Satz~2 of \S 19 in~\cite{Fr55}. 

\proclaim{Theorem 13}
Let $J$ be an invariant almost complex structure on $G/H$ and 
$\Vc$ a closed linear subspace of $\gg$ such that 
$\gg=\hg\dotplus\Vc$. 
Now let $I\in\Ic_{\Vc}(G,H)$ $(\subseteq\Bc(\gg))$ with $c_{\Vc}(I)=J_{p_0}$. 
Then the following assertions are equivalent: 

\itemitem{\rm(i)} There exists on $G/H$ a structure of complex 
Banach manifold whose underlying real analytic structure is 
just the original structure of $G/H$ and such that 
$J_pv=\ir v$ whenever $p\in G/H$ and $v\in T_p(G/H)$.

\itemitem{\rm(ii)} For all $x,y\in\gg$ we have 
$I[Ix,y]+I[x,Iy]+[x,y]-[Ix,Iy]\in\hg$. 

\endproclaim

\demo{Proof}
We begin with a few general remarks. 
First note that the existence of $I$ in the statement follows by 
Propositions 10~and~12. 
Next, for $x\in\gg$ arbitrary, we introduce the notation 
$\widetilde{x}$ for the corresponding $G$-invariant (real analytic) vector field 
on $G/H$. 
Note that for each $x\in\gg$ we have 
$$T_{\1}\pi(x)=\widetilde{x}_{p_0}\text{ and }\widetilde{Ix}=J\widetilde{x}.
\eqno{(3)}$$
Moreover, since $I\in\Ic_{\Vc}(G,H)$ (see Proposition~12), we have 
$T_{\1}\pi\circ I\circ\Ad_G(h)|_{\Vc}=T_{\1}\pi\circ\Ad_G(h)\circ I|_{\Vc}$ 
for all $h\in H$. 
By differentiating this equality with respect to $h\in H$, we see that 
$$T_{\1}\pi(I[z,v])=T_{\1}\pi[z,Iv]\text{ whenever }z\in\hg\text{ and }v\in\Vc.$$
Also, recall from the commutative diagram in the proof of Proposition~12 that 
$T_{\1}\pi\circ I=J_{p_0}\circ T_{\1}\pi$. 
This implies that for all $x,y\in\Vc$ we have 
$$\aligned
T_{\1}\pi(I[Ix,y]+&I[x,Iy]+[x,y]-[Ix,Iy]) \cr
 &=J_{p_0}(T_{\1}\pi[Ix,y])+J_{p_0}(T_{\1}\pi[x,Iy])
   +[\widetilde{x},\widetilde{y}]_{p_0}-\widetilde{[Ix,Iy]}_{p_0} \cr
 &=J_{p_0}\widetilde{[Ix,y]}_{p_0}+J_{p_0}\widetilde{[x,Iy]}_{p_0}
   +[\widetilde{x},\widetilde{y}]_{p_0}-[\widetilde{Ix},\widetilde{Iy}]_{p_0} \cr
 &=J_{p_0}[\widetilde{Ix},\widetilde{y}]_{p_0}+J_{p_0}[\widetilde{x},\widetilde{Iy}]_{p_0}
   +[\widetilde{x},\widetilde{y}]_{p_0}-[J\widetilde{x},J\widetilde{y}]_{p_0} \cr
 &=(J[J\widetilde{x},\widetilde{y}]+J[\widetilde{x},\widetilde{Jy}]
   +[\widetilde{x},\widetilde{y}]-[J\widetilde{x},J\widetilde{y}])_{p_0},
\endaligned$$
where the last equality follows by (3). 

We now proceed to prove that assertions (i)~and~(ii) are equivalent. 

(i)$\Rightarrow$(ii) 
It follows by Theorem~7 that the torsion of $J$ vanishes, 
hence the above computation shows that for all $x,y\in\Vc$ we have 
$T_{\1}\pi(I[Ix,y]+I[x,Iy]+[x,y]-[Ix,Iy])=0$, 
that is, $I[Ix,y]+I[x,Iy]+[x,y]-[Ix,Iy]\in\hg$. 

Moreover, if $x\in\hg$ and $y\in\Vc$, then 
$$\aligned
T_{\1}\pi(I[Ix,y]+&I[x,Iy]+[x,y]-[Ix,Iy]) \cr
 &=T_{\1}\pi[Ix,Iy]+T_{\1}\pi[x,I^2y]+T_{\1}\pi[x,y]-T_{\1}\pi[Ix,Iy] \cr
 &=0,
\endaligned$$
hence also $I[Ix,y]+I[x,Iy]+[x,y]-[Ix,Iy]\in\hg$. 

Finally, if $x,y\in\hg$ then 
$I[Ix,y]+I[x,Iy]+[x,y]-[Ix,Iy]=[x,y]\in\hg$.

(ii)$\Rightarrow$(i) 
By the computation performed just before the proof of (i)$\Rightarrow$(ii),  
for $x,y\in\Vc$ arbitrary we get 
$(J[J\widetilde{x},\widetilde{y}]+J[\widetilde{x},\widetilde{Jy}]
   +[\widetilde{x},\widetilde{y}]-[J\widetilde{x},J\widetilde{y}])_{p_0}=0$. 
Since the almost complex structure $J$ is invariant and both 
$\widetilde{x}$ and $\widetilde{y}$ are $G$-invariant vector fields on $G/H$, 
it then follows that 
$J[J\widetilde{x},\widetilde{y}]+J[\widetilde{x},\widetilde{Jy}]
   +[\widetilde{x},\widetilde{y}]-[J\widetilde{x},J\widetilde{y}]=0$ on $G/H$. 
But this implies that the torsion of $J$ vanishes, since 
for every $p\in G/H$ we have 
$T_p(G/H)=\{\widetilde{x}\mid x\in\Vc\}$ in view of the fact 
$\pi\colon G\to G/H$ is a submersion and 
$\Vc\dotplus\Ker(T_{\1}\pi)=\gg$. 
\qed
\enddemo

\definition{Definition 14}
We denote by $\Ico(G,H)$ the set of all invariant almost complex structures $J$ 
on $G/H$ such that assertion~(i) in Theorem~13 holds. 
Each element $J\in\Ico(G,H)$ will be called an {\it invariant complex structure} 
on the homogeneous space $G/H$. 
\enddefinition

Now we are able to give the following parameterization of the invariant complex structures in Lie algebraic terms. 
When specialized to finite dimensions, 
the next theorem becomes a more precise version of 
the Satz of \S 20 in~\cite{Fr55}. 

\proclaim{Theorem 15}
Denote by $\Kc_0(G,H)$ the set of all closed complex subalgebras $\kg$ of 
$\gg_{\Bbb C}$ 
such that 
$$\kg+\overline{\kg}=\gg_{\Bbb C},\;
\kg\cap\overline{\kg}=\hg_{\Bbb C}
\text{ and }
\Ad_G(h)\kg\subseteq\kg
\text{ for all }h\in H.$$
Next fix a closed linear subspace $\Vc$ of $\gg$ such that $\gg=\hg\dotplus\Vc$ 
and for each $J\in\Ico(G,H)$ 
define 
$$b_{\Vc}(J)=\hg_{\Bbb C}+\Ker(I-\ir\cdot\id_{\Vc_{\Bbb C}}),$$
where $I:=(c_{\Vc})^{-1}(J_{p_0})\in\Ic_{\Vc}(G,H)\subseteq\Bc(\gg)$. 

Then we obtain a bijection
$$b_{\Vc}\colon\Ico(G,H)\to\Kc_0(G,H).$$
\endproclaim

\demo{Proof}
The proof has several stages. 

$1^\circ$ 
At this stage we show that $b_{\Vc}$ is well defined in the sense that 
for $J\in\Ico(G,H)$ arbitrary we have indeed 
$\kg:=c_{\Vc}(J)\in\Kc_0(G,H)$. 
What is obvious is that $\kg$ is a closed complex linear subspace of $\gg_{\Bbb C}$. 

Since $I^2=-\id_{\Vc_{\Bbb C}}$ on $\Vc_{\Bbb C}$, 
we get $\Vc_{\Bbb C}
=\Ker(I-\ir\cdot\id_{\Vc_{\Bbb C}})\dotplus\Ker(I+\ir\cdot\id_{\Vc_{\Bbb C}})$, 
whence 
$$\gg_{\Bbb C}=\hg_{\Bbb C}\dotplus\Vc_{\Bbb C}
=\hg_{\Bbb C}
\dotplus\Ker(I-\ir\cdot\id_{\Vc_{\Bbb C}})\dotplus\Ker(I+\ir\cdot\id_{\Vc_{\Bbb C}}). 
\eqno{(4)}$$
Now, since $I\overline{z}=\overline{Iz}$ for all $z\in\gg_{\Bbb C}$ 
and $\kg=\hg_{\Bbb C}+\Ker(I-\ir\cdot\id_{\Vc_{\Bbb C}})$, 
it follows by~(4) that 
$\kg+\overline{\kg}=\gg_{\Bbb C}$ and 
$\kg\cap\overline{\kg}=\hg_{\Bbb C}$. 

Furthermore let $h\in H$ arbitrary. 
If $v\in\hg_{\Bbb C}$ it is clear that $\Ad_g(h)v\in\hg_{\Bbb C}$. 
Now let $v\in\Ker(I-\ir\cdot\id_{\Vc})$ arbitrary. 
Since $I\in\Ic_{\Vc}(G,H)$, 
it follows that $(I\circ\Ad_G(h))v-(\Ad_g(h)\circ I)v\in\hg_{\Bbb C}$, 
that is, 
$I(\Ad_G(h)v)-\ir\cdot\Ad_G(h)v\in\Ker I$. 
Consequently
$0=I(I(\Ad_G(h)v)-\ir\cdot\Ad_G(h)v)
=-\Ad_G(h)v-\ir\cdot(I\circ\Ad_G(h))v$, 
whence $I(\Ad_G(h)v)=\ir\cdot\Ad_G(h)v$, 
which shows that $\Ad_G(h)v\in\Ker(I-\ir\cdot\id_{\Vc})$. 
Thus $\Ad_G(h)\kg\subseteq\kg$ for arbitrary $h\in H$. 

Thence by differentiating with respect to $h\in H$ we get 
$[\hg,\kg]\subseteq\kg$, whence $[\hg_{\Bbb C},\kg]\subseteq\kg$. 
Thus, in order to prove that $\kg$ is a subalgebra of $\gg_{\Bbb C}$, 
it remains to check that $[x,y]\in\kg$ for arbitrary 
$x,y\in\Ker(I-\ir\cdot\id_{\Vc})$. 
In fact, since $J\in\Ico(G,H)$, it follows by Theorem~13((i)$\Rightarrow$(ii)) 
that $I[Ix,y]+I[x,Iy]+[x,y]-[Ix,Iy]\in\hg_{\Bbb C}$, 
that is, 
$I[\ir\cdot x,y]+I[x,\ir\cdot y]+[x,y]-[\ir\cdot x,\ir\cdot y]\in\hg_{\Bbb C}$, 
whence $2\ir\cdot I[x,y]+2[x,y]\in\hg_{\Bbb C}=\Ker I$. 
This implies that 
$I(2\ir\cdot I[x,y]+2[x,y])=0$, 
whence $-2\ir[x,y]+2I[x,y]=0$, and thus 
$[x,y]\in\Ker(I-\ir\cdot\id_{\Vc_{\Bbb C}})$.

$2^\circ$ 
At this stage we just remark that, for arbitrary $\kg\in\Kc_0(G,H)$, 
the natural embedding $\gg\hookrightarrow\gg_{\Bbb C}$ 
induces an isomorphism of Banach spaces 
$$\nu_{\kg}\colon\gg/\hg\to\gg_{\Bbb C}/\kg,\quad x+\hg\mapsto x+\kg.$$
In fact, since 
$\kg\cap\gg=\kg\cap\overline{\kg}\cap\gg=\hg_{\Bbb C}\cap\gg=\hg$, 
it follows that $\Ker\nu_{\kg}=\{0\}$. 
On the other hand, to prove that $\nu_{\kg}$ is surjective, 
let $z\in\gg_{\Bbb C}$ arbitrary. 
Then there exist $y_1,y_2\in\kg$ such that $x=y_1+\overline{y_2}$, 
whence $z=(y_2+\overline{y_2})+(y_1-y_2)\in\gg+\kg$, 
and this shows that $z+\kg=(y_2+\overline{y_2})+\kg\in\Ran\nu_{\kg}$. 

$3^\circ$ 
We now use the remark of stage~$2^\circ$ to define a mapping 
$$\beta\colon\Kc_0(G,H)\to\Ic_0(G,H)$$
such that for each $\kg\in\Kc_0(G,H)$ the diagram 
$$\CD
\gg_{\Bbb C}/\kg @>{\mult}>> \gg_{\Bbb C}/\kg \cr
@A{\nu_{\kg}}AA @AA{\nu_{\kg}}A \cr 
\gg/\hg @>{\beta(\kg)}>> \gg/\hg
\endCD$$
is commutative, where $\mult$ stands for the multiplication-by-i operator 
on the complex vector space $\gg_{\Bbb C}/\kg$. 

It is not hard to see that the mapping $\beta$ indeed takes 
values in $\Ic_0(G,H)$. 
In fact, let $\kg\in\Kc_0(G,H)$ arbitrary. 
We have $(\mult)^2=-\id_{\gg_{\Bbb C}/\kg}$, 
hence $(\beta(\kg))^2=-\id_{\gg/\hg}$. 
Moreover, for each $h\in H$ we have 
$\Ad_G(h)\kg\subseteq\kg$, 
hence $\Ad_G(h)(\ir\cdot z+\kg)=\ir\cdot\Ad_G(h)z+\kg$ 
for all $z\in\gg_{\Bbb C}$. 
Then $\Adt_G(h)\circ\beta(\kg)=\beta(\kg)\circ\Adt_G(h)$ for all $h\in H$, 
and this concludes the proof that $\beta(\kg)\in\Ic_0(G,H)$. 
(We recall from the first remark in stage~$1^\circ$ of 
the proof of Proposition~12 that $\Adt_G(h)=T_{p_0}(\alpha_h)$ 
for $h\in H$.) 

$4^\circ$ 
Let $\kg\in\Kc_0(G,H)$ arbitrary. 
We are going to prove that 
$$J:=\widehat{\beta(\kg)}\in\Ico(G,H),$$
which will eventually show that the mapping 
$\widehat{\beta(\cdot)}\colon\Kc_0(G,H)\to\Ico(G,H)$ 
is an inverse to $b_{\Vc}$. 

We have seen in stage~$3^\circ$ that $\beta(\kg)\in\Ic_0(G,H)$, 
hence Proposition~10 shows that $J\in\Iac(G,H)$. 
Moreover, 
we have $\gg_{\Bbb C}=\hg_{\Bbb C}\dotplus\Vc_{\Bbb C}$ and 
$\hg_{\Bbb C}=\kg\cap\overline{\kg}$, 
hence $\kg=\hg\dotplus(\Vc_{\Bbb C}\cap\kg)$ and 
$\overline{\kg}=\hg\dotplus(\Vc_{\Bbb C}\cap\overline{\kg})$. 
Since $\gg_{\Bbb C}=\kg+\overline{\kg}$, it follows that 
$$\gg_{\Bbb C}
=\hg_{\Bbb C}\dotplus
(\Vc_{\Bbb C}\cap\kg)\dotplus(\Vc_{\Bbb C}\cap\overline{\kg}).$$
Using this decomposition of $\gg_{\Bbb C}$, 
now define a ${\Bbb C}$-linear mapping $I\colon\gg_{\Bbb C}\to\gg_{\Bbb C}$ 
such that $\hg_{\Bbb C}=\Ker I$, 
$\Vc_{\Bbb C}\cap\kg=\Ker(I+\ir\cdot\id_{\gg_{\Bbb C}})$ and 
$\Vc_{\Bbb C}\cap\overline{\kg}=\Ker(I-\ir\cdot\id_{\gg_{\Bbb C}})$. 

In particular we have $\overline{Iz}=I\overline{z}$ for all $z\in\gg_{\Bbb C}$, 
hence $I(\gg)\subseteq\gg$. 
We shall check shortly that $I|_{\gg}\in\Ic_{\Vc}(G,H)$. 
In fact $(I|_{\Vc})^2=-\id_{\Vc}$ since 
$(I|_{\Vc_{\Bbb C}})^2=-\id_{\Vc_{\Bbb C}}$. 
Next note that for $h\in H$ and $z\in\Vc\cap\kg\subseteq\kg$ 
we have $\Ad_G(h)z\in\kg=\hg\dotplus(\Vc_{\Bbb C}\cap\kg)$ hence 
there exist $z_{\hg_{\Bbb C}}\in\hg_{\Bbb C}$ and $z_{\Vc_{\Bbb C}}\in\Vc_{\Bbb C}$ 
such that $\Ad_G(h)z=z_{\hg_{\Bbb C}}+z_{\Vc_{\Bbb C}}$, 
whence 
$$I(\Ad_G(h)z)=I(z_{\hg_{\Bbb C}}+z_{\Vc_{\Bbb C}})
=\ir\cdot z_{\Vc_{\Bbb C}}.$$
Since $\Ad_G(h)Iz=\Ad_G(h)(\ir\cdot z)=\ir\cdot\Ad_G(h)z
=\ir\cdot z_{\hg_{\Bbb C}}+\ir\cdot z_{\Vc_{\Bbb C}}$, 
it follows that 
$I(\Ad_G(h)z)-\Ad_G(h)Iz\in\hg_{\Bbb C}$. 
In particular
$T_{\1}\pi\circ I\circ\Ad_G(h)|_{\Vc}=T_{\1}\pi\circ\Ad_G(h)\circ I|_{\Vc}$, 
and thus 
$$I\in\Ic_{\Vc}(G,H).$$
Also, 
$$c_{\Vc}(I)=\beta(\kg)=J_{p_0}.$$
Thus, in order to show that $J\in\Ico(G,H)$, we have to show that $I$ satisfies 
condition~(ii) in Theorem~13. 
To this end, it is convenient to make use of the natural projection 
$$\chi\colon\gg_{\Bbb C}=\Vc_{\Bbb C}\dotplus\hg_{\Bbb C}\to\hg_{\Bbb C}$$
with $\Ker\chi=\Vc_{\Bbb C}$, since we have  
$$
Ix=\cases 
 I(x-\chi(x))=-\ir(x-\chi(x))=-\ir\cdot x+\ir\cdot\chi(x)
\text{ if }x\in\kg
%\eqno{(5)} 
\cr 
 I(x-\chi(x))=\ir(x-\chi(x))=\ir\cdot x-\ir\cdot\chi(x)
\text{ if }x\in\overline{\kg}.
%\eqno{(6)}
\endcases$$
Now, using this fact, it is straightforward do show that 
$I[Ix,y]+I[x,Iy]+[x,y]-[Ix,Iy]\in\hg_{\Bbb C}$ 
for all $x,y\in\gg_{\Bbb C}$, 
whence 
$I[Ix,y]+I[x,Iy]+[x,y]-[Ix,Iy]\in\hg$ 
for $x,y\in\gg$. 
Thus $I$ satisfies condition~(ii) in Theorem~13, 
and then $J\in\Ico(G,H)$. 

$5^\circ$ 
Now consider the mapping 
$\widehat{\beta(\cdot)}\colon\Kc_0(G,H)\to\Ico(G,H)$, 
which is well defined according to stage~$4^\circ$. 
It is easy to see that this mapping is inverse to 
$b_{\Vc}\colon\Ico(G,H)\to\Kc_0(G,H)$, 
hence $b_{\Vc}$ is a bijection and the proof ends. 
\qed
\enddemo

\proclaim{Corollary 16 {\rm(\cite{MS98})}}
Let $A$ be a unital $C^*$-algebra, $n\ge1$ an integer, and denote by $\Fl_{A,n}$ the set of all $n-tuples$ $(p_1,\dots,p_n)$ 
of orthogonal projections in $A$ satisfying $p_1\le\cdots\le p_n$. 
Then $\Fl_{A,n}$ has a structure of complex manifold such that for every unitary element $u\in A$ the mapping 
$$\Fl_{A,n}\to\Fl_{A,n},\quad 
(p_1,\dots,p_n)\mapsto(up_1u^*,\dots,up_nu^*),$$
is holomorphic. 
\endproclaim

\demo{Proof}
Denote by $\U_A$ the unitary group of $A$ 
and for each $(p_1,\dots,p_n)\in\Fl_{A,n}$ define 
$$\Fl_A(p_1,\dots,p_n)=\{(up_1u^*,\dots,up_nu^*)\mid u\in\U_A\}.$$
Clearly it suffices to equip each set $\Fl_A(p_1,\dots,p_n)$ 
with a $\U_A$-invariant complex structure, 
since $\Fl_{A,n}$ is the disjoint union of the sets of this form. 

Next fix $(p_1,\dots,p_n)\in\Fl_{A,n}$ and denote $p_0=0$ and 
$$\aligned
G&=\U_A,\quad\gg=\{a\in A\mid a^*=-1\},\cr
H&=\{u\in\U_A\mid up_j=p_ju\text{ for }1\le j\le n\},\quad 
\hg=\{a\in\gg\mid up_j=p_ju\text{ for }1\le j\le n\}.
\endaligned$$
Then it is clear that $H$ is a Banach-Lie subgroup of $G$ 
(since it is closed in $G$, the exponential map of $G$ gives by restriction 
a homeomorphism of a neighborhood of $0\in\hg$ onto a neighborhood of 
$\1\in H$, 
and $\{a\in\gg\mid p_jap_l=0\text{ for }j\ne l\}$ is a direct complement to $\hg$ in $\gg$). 
Moreover, denoting by 
$$b\mapsto \overline{b}=-b^*$$
the conjugation of $\gg_{\Bbb C}=A$ whose fixed-point set is $\gg$, 
it follows at once that 
$$\kg:=\{b\in A\mid (p_i-p_{i-1})b(p_j-p_{j-1})=0\text{ if }
1\le j<i\le n\}$$
is an element of $\Kc_0(G,H)$, with the notation of Theorem~15. 
Then Theorem~15 implies that there exists on 
$G/H\simeq\Fl_A(p_1,\dots,p_n)$ a $G$-invariant complex structure. 
Since $G=\U_A$, the proof is finished. 
\qed
\enddemo

\bigskip

{\it Acknowledgment.} 
We are indebted to Professor L\'aszl\'o Lempert for 
a lot of useful discussions and in particular for drawing our attention 
to the paper~\cite{Pa00}. 
We also thank Professor Tudor Ratiu for several illuminating comments on the Frobenius theorem 
in infinite dimensions. 

\head References\endhead

\refstyle{ABCDE}

\enddocument